\documentclass[twoside,11pt]{article}

\usepackage{amsfonts,amsmath,amssymb,enumerate}
\usepackage{euscript,mathrsfs}
\usepackage{color,graphics}
\usepackage{empheq}

\numberwithin{equation}{section}

\usepackage{titlesec}
\titleformat{\section}{\bfseries}{\thesection}{1em}{}
\titleformat{\subsection}{\itshape}{\thesubsection}{1em}{}

\newfont{\ctv}{msam10}
\newcommand{\bbox}{\mbox{\ctv \symbol{4}}}
\def\QED{{$\hfill\bbox$}}
\newenvironment{pf}[1]{\par\vskip1mm{\noindent\it #1.}\ }{\QED\par\vskip2mm}

\def\dd{\,\mathrm{d}}

\definecolor{vio}{rgb}{0.7,0,0.8}

\newtheorem{theorem}{Theorem}[section]
\newtheorem{lemma}[theorem]{Lemma}
\newtheorem{corollary}[theorem]{Corollary}
\newtheorem{proposition}[theorem]{Proposition}
\newtheorem{remark}[theorem]{Remark}

\numberwithin{equation}{section}

\def\bproof{\noindent{\it Proof. }}
\def\eproof{\null\hfill {$\Box$}\bigskip}

\def\bbproof{\noindent{\it End of the proof of  }}
\def\eeproof{\null\hfill {$\Box$}\bigskip}

\def\Vo{\vbox{\offinterlineskip\hbox{\kern 3pt$\scriptstyle\circ$}
\kern 1pt\hbox{$V$}}}

\def\Ho{\vbox{\offinterlineskip\hbox{\kern 3pt$\scriptstyle\circ$}
\kern 1pt\hbox{$H$}}}

\def\Wo{\vbox{\offinterlineskip\hbox{\kern 3pt$\scriptstyle\circ$}
\kern 1pt\hbox{$W$}}}

\newcommand{\bt}{\begin{theorem}}
\newcommand{\et}{\end{theorem}}

\def\be{\begin{equation}\label}
\def\ee{\end{equation}}

\def\barr{\begin{array}}
\def\earr{\end{array}}

\def\ber{\begin{eqnarray}}
\def\eer{\end{eqnarray}}

\def\bers{\begin{eqnarray*}}
\def\eers{\end{eqnarray*}}

\def\bpf{\begin{pf}}
\def\epf{\end{pf}}

\def\eqldef{\overset{\text{\tiny \rm def}}{=}}

\newcommand{\abs}[2][{}]{\lvert#2\rvert_{#1}}
\def\dd{\;\!\mathrm{d}}    

\newcommand{\Er}{\mathbb{R}}
\newcommand{\En}{\mathbb{N}}

\newcommand{\norm}[2][{}]{\lVert#2\rVert_{{#1}}}

\renewcommand{\tilde}{\widehat}
\renewcommand{\tilde}{\widetilde}

\def\eqldef{\overset{\text{\tiny\rm def}}{=}}
\def\L{{\rm{L}}}
\def\H{{\rm{H}}}

\def\E{{\rm{E}}}

\def\C{{\mathrm{C}}}
\def\H{{\mathrm{H}}}
\def\L{{\mathrm{L}}}

\begin{document}

\title{Weak solutions for a singular beam equation}

\author{Olena~Atlasiuk\thanks{Institute of Mathematics Czech Academy of Sciences \& 
Institute of Mathematics of the National Academy of Sciences of Ukraine,
\v{Z}itn\'{a} 25, CZ-11567 Praha 1, Czech Republic,
(\tt  atlasiuk@math.cas.cz)}\and Arnaud~Heibig\thanks{Univ Lyon, INSA Lyon, UJM, UCBL, ECL, CNRS UMR 5208, ICJ, 
43 boulevard du 11 novembre 1918, 
F--69622 Villeurbanne Cedex, France 
(\tt  arnaud.heibig@insa-lyon.fr, adrien.petrov@insa-lyon.fr)}
\and Adrien~Petrov$^{\dagger}$}

\pagestyle{myheadings}
\thispagestyle{plain}
\markboth{O.~Atlasiuk, A.~Heibig, A.~Petrov}
{\small{Weak solutions for a singular beam equation}}


\maketitle

\vspace{-2em}

\begin{abstract}
\noindent This paper deals with a dynamic Gao beam of infinite length subjected to a moving concentrated Dirac mass.
Under appropriate regularity assumptions on the initial data, the problem possesses a weak solution which is obtained 
as the limit of a sequence of solutions of regularized problems.
    
\end{abstract}


\hspace*{-0.6cm}{\textbf {Key words.}} Gao beam equation, Dirac mass,  moving load, existence result, Sobolev spaces.

\hspace*{-0.6cm}{\textbf{AMS Subject Classification:}} 35B45, 35Q74, 74H20, 74K10, 74M15
 
\vspace{0.3cm}


\section{Description of the problem}
\label{descritption}

The behavior of a beam plays a crucial role in various applications as in railway track design. 
The railway compagnies aim to enhance the train rolling speed to meet the increased demands 
in both passenger and freight transportation worldwide. Identifying the factors contributing to the 
occurrence of railway track defects is rather important to maintain the required track quality. The
oscillations amplitudes in railway tracks due to the trains movement is intensively studied in 
scientific literature (see \cite{Fry72VSSM,DimRod12CVUM}). These oscillations lead to undesirable 
consequences such as a premature wear and deformation of railway tracks. Understanding the impact 
of these oscillations on system reliability is essential to maintain the track quality necessary for traffic safety 
and passenger comfort.

Most of mathematical models consider rails as beams under the influence 
of moving loads and their properties are intensively studied. A considerable engineering and mathematical 
literature is devoted to study to numerous mathematical models for beams, such as the  
Euler-Bernoulli (see \cite{DekNic99CEBB,CasZua00EBEBI,AhnSte05EBDC,HeKrPe22SSEB}), Timoshenko (see \cite{KimRen87BCTB,BeLoRe12ATBI}) and Gao beams. 
The Euler-Bernoulli and Timoshenko beam theories are the oldest ones and they are widely employed nowadays
in various structural analysis methodologies in engineering. The Euler-Bernoulli beam is commonly used to model the 
bending beams behavior, the axial and shear deformations are neglected compared to bending.
Thus, the beam cross-section remains planar during loading.  
The Timoshenko beam is suitable to study the shear beams behavior, thereby providing 
an accurate representation of deformation within the cross-section of the beam.

The current study involves analysis of the horizontal motion of a vertical point-load on a metallic rail, considering 
scenarios where the load has mass or is massless, and examining the resulting oscillations of the system. 
The mathematical model of the Gao beam was originally introduced in \cite{Gao96NBTC}. 
However the Gao beam is also studied in different contexts (see \cite{AMS09VNDB,AhKuSh12DCGB,AnLiSh15EDCS,MacNet15SCPG,GaMaNe15MFNG,MacNet18SCPE}). 
Notice also that semilinear beam is studied in \cite{TakYos12SBAP,TakYos13SBWD}.
Comparisons between simulations of the Gao beam and the Euler-Bernoulli linear beam reveal significant 
differences, indicating that the linear beam is suitable only for small loads, whereas the Gao beam accommodates moderate loads (see \cite{DBK*19VGBM}).

We consider in this paper a straight infinite Gao beam of thickness \(2h\). 
A detailed derivation of the Gao beam model can be found in \cite{Gao96NBTC,AMS09VNDB}. An horizontal traction \(p\), which is a time-depending function, 
is applied at one end. We are focused in this work on examining the sideways movement of a concentrated load \(P(t)\), where the load may 
change over time. This load is applied to a mobile mass \(m\) positioned at \(\zeta(t)\) with a horizontal velocity 
\(\dot\zeta(t)\), induced by a horizontal applied force. 
The transverse displacement of the Gao beam \(u(t,x)\) for \((t,x)\in(0,T)\times\Er\) is governed by a following partial differential equation: 
\begin{equation}
\label{AtHePe:mod}
\begin{aligned}
&\varrho u_{tt}+m\delta(x-\zeta(t))u_{tt}-m\delta'(x-\zeta(t))\dot\zeta(t)u_t+ku_{xxxx}
-(eu_x^2-\nu p)u_{xx}\\&=\varrho f+\delta(x-\zeta(t))P(t),
\end{aligned}
\end{equation}
where \(\delta(\cdot)\) is the Dirac function,
\(\varrho\) is the material density,
\(f\) is the applied mechanical loading, \(k\eqldef \frac{2h^3 E}{3(1-\bar\nu^2)}\), 
\(\nu\eqldef (1+\bar\nu)\) and \(e\eqldef 3hE\) where \(E\) and \(\bar\nu\) denote the Young modulus and
the Poisson ratio, respectively. 
Here and below 
\((\cdot)_t\eqldef \partial_t(\cdot)\),
\((\cdot)_x\eqldef\partial_x(\cdot)\)
denote the partial derivatives with respect to \(t\) and \(x\), respectively, while
\(\dot{(\cdot)}\) and \((\cdot)'\) denote the derivatives with respect to \(t\) and \(x\), respectively.
We prescribe also initial data
\begin{equation}
\label{AtHePe:Init_cond}
u(0,\cdot)=u_0\quad\text{and}\quad u_t(0,\cdot)=u_1.
\end{equation}
For further explanations on this model, the reader is referred to \cite{DBK*19VGBM} and to the references therein.
Equation \eqref{AtHePe:mod} can be considered the limit of the equation treated in \cite{DBK*19VGBM} when the gamma viscosity tends to zero, and is related 
to the usual viscous beam equation without Dirac measure in the time derivative. Nevertheless, this leads to an entirely different 
mathematical problem since three-order estimates cannot be obtained using the viscous term, and what is more, 
since the usual energy estimates fail. 
In particular, the non-linear term should be addressed
 somewhat indirectly since no \(\L^{\infty}\) estimate is available immediately.

The positive constants \(\varrho\), \(m\),
\(k\) and \(e\) play any role in the mathematical analysis carried out below. Consequently, without loss of generality,
we set them equal to \(1\). Notice, however, that the case of non constant coefficients
(for example \(k=k(x)\)) could be processed using approximate square roots for elliptic operators.

Throughout this paper, we assume that \(\zeta\in\C^2([0,T])\), \(P\in\C^2([0,T])\), \(p\in\C^0([0,T];\H^2(\Er))\) and \(f\in\C^0([0,T];\H^2(\Er))\), but these assumptions
could be weakened.
The main result of this paper can be stated as follows.

\begin{theorem}
\label{AtHePe:thm1}
Let \(T>0\). Assume that \(u_0\in\H^2(\Er)\), \(u_1\in\H^1(\Er)\), \(\zeta\in\C^2([0,T])\), \(P\in\C^2([0,T])\) and \(f,p\in\C^0([0,T];\H^2(\Er))\). 
Then there exists a function 
\(u\in\C^0([0,T);\H_{\rm{loc}}^{2-\alpha}(\Er))\cap\L^{\infty}([0,T];\H^2(\Er))\), for any \(\alpha\in\,]0,2[\,\), and \(u_t\in\L^{\infty}([0,T];\L^2(\Er))\) 
such that for any \(v\in\C^{2}([0,T];\\\H^2(\Er))\) with a compact support in \([0,T[\,\times\Er\), we have
\begin{equation}
\label{AtHePe:FV}
\begin{cases}
\displaystyle{\int_0^T\int_{\Er} u_t v_t\dd x\dd t+\int_0^T\int_{\Er} u_{xx} v_{xx}\dd x\dd t
+\frac13\int_0^T\int_{\Er}u_x^3 v_x\dd x\dd t}
\\\noalign{\vspace*{\jot}} 
\displaystyle{-\int_0^T\int_{\Er}\nu p u_{xx} v \dd x\dd t-u_0(\zeta(0))v_t(0,\zeta(0))
-\int_0^T u(t,\zeta(t))v_{tt}(t,\zeta(t))\dd t}
\\\noalign{\vspace*{\jot}} 
\displaystyle{-\int_0^T\dot\zeta(t)u_x(t,\zeta(t))v_t(t,\zeta(t))\dd t
-\int_0^T\dot\zeta(t)u(t,\zeta(t))v_{xt}(t,\zeta(t))\dd t} 
\\\noalign{\vspace*{\jot}} 
\displaystyle{+\int_0^T\int_{\Er} f v\dd x\dd t
+\int_0^TP(t)v(t,\zeta(t))\dd t
+\int_{\Er}u_1 v(0,\cdot)\dd x}\\
+u_1(\zeta(0))v(0,\zeta(0))=0.
\end{cases}
\end{equation}
Moreover, \(u(0,\cdot)=u_0\) holds.
\end{theorem}
In this statement, we refrain from expressing the term
\(\int_0^T\int_{\Er}\frac{\dd}{\dd t}(\delta(x-\xi(t))\partial_t u(t,\xi(t)))v(t,x)\dd x\dd t\) as 
\(-\int_0^T \partial_tu(t,\xi(t))\partial_tv(t,\xi(t))\dd t\). Indeed, the trace is not defined within our functional frame. 
To establish such an expression, higher-order regularity would be required. 
However, in the linear homogenous case, 
the formula \(\partial_t^{(2)}u(t,x)+\partial_x^{(4)}u(t,x)=-\frac{\rm{d}}{\rm{d}t}(\delta(x-\zeta(t))\partial_t u(t,x))\) entails that \(\partial_t^{(2)}u\) and \(\partial_x^{(4)}u\)
cannot both possess \(\L^2(0,T;\L^2(\Er))\) regularity. The same formula suggests that 
\(\partial_t^{(2)}u\notin \C^0([0,T];\L^2(\Er))\)
as the formal initial data
\(\partial_t^{(2)}u(0,x)=-(\partial_x^{(4)}u(0,x)-\frac{\rm{d}}{\rm{d}t}(\delta(x-\zeta(0))\partial_tu(0,x)))\)
does not belong to \(\L^2(\Er)\) for a smooth data \(u_0(x)=u(0,x)\). 
Lastly, the usual energy estimate, obtained by multiplying \eqref{AtHePe:mod} by \(\partial_t u\) and by integrating 
the resulting result over \([0,T]\times\Er\), can not be performed due to  the term \(\delta(x-\zeta(t))\). 

The paper is organized as follows. An existence result for a linear mollified problem is presented in Section \ref{prelim}.
Section \ref{Approx_Problem} addresses a mollified  nonlinear equation; some uniform a priori estimates are derived. The proof of Theorem \ref{AtHePe:thm1} is provided 
in Section \ref{Passage_limit}.


\section{A preliminary existence result}
\label{prelim}


Under some regularity assumptions on the data, we present an existence and uniqueness result
to the system \eqref{AtHePe:2} in a suitable function space. To this aims,
for any \(t>0\), we introduce the following sets:
\begin{subequations}
\label{AtHePe:1}
\begin{align}
\mathrm{E}_2(t)&\eqldef \C^0([0,t];\H^2(\Er)),\\
\mathrm{E}_1(t)&\eqldef \C^0([0,t];\H^2(\Er))\cap \C^1([0,t];\L^2(\Er)),\\
\mathrm{E}_0 (t)&\eqldef \C^0([0,t];\H^4 (\Er))\cap \C^1([0,t];\H^2 (\Er))\cap \C^2([0,t];\L^2 (\Er)).
\end{align}
\end{subequations}
Hence, we have \(\mathrm{E}_0(t)\hookrightarrow\mathrm{E}_1(t)\hookrightarrow\mathrm{E}_2(t)\).
Next, let \(\theta\in\mathcal{D}(\Er)\) be an even density of probability. 
For any \((t,x)\in\Er_{+}\times\Er\), we set
\begin{equation*}
F(t,x)\eqldef \frac1{1+\theta(x-\zeta (t))}\quad\text{and}\quad
G(t,x)\eqldef -\frac{\dot\zeta(t)\theta'(x-\zeta(t))}{1+\theta(x-\zeta(t))}.
\end{equation*}
We denote below by \(C_{\eta}>0\) a generic constant depending on \(\eta\). 
\begin{proposition}
\label{AtHePe:Prop1}
Let \(T>0\), \(u_0\in\H^4(\Er)\), \(u_1\in\H^2(\Er)\) and \(g\in \mathrm{E}_2(T)\). Then, there 
exists a unique \(u\in\mathrm{E}_0(T)\) solution to the following system:
\begin{subequations}
\label{AtHePe:2}
\begin{empheq}[left=\empheqlbrace]{alignat=1}
\label{AtHePe:2_a}
&u_{tt}+Fu_{xxxx}+G u_t=g,\\&
\label{AtHePe:2_b}
u(0,\cdot)=u_0\quad\text{and}\quad u_t(0,\cdot)=u_1.
\end{empheq}
\end{subequations}
Moreover, \(u\) satisfies the following inequality:
\begin{equation}
\label{AtHePe:3}
\norm[\H^2(\Er)]{u(t,\cdot)}+
\norm[\L^2(\Er)]{u_t(t,\cdot)}
\leq
C_T\bigl(\norm[\H^2(\Omega)]{u_0}+\norm[\L^2(\Omega)]{u_1}+\norm[\L^2(0,t;\L^2(\Er))]{g}\bigr).
\end{equation}
\end{proposition}
The proof is omitted. It follows from Proposition \ref{AtHePe:Prop1} that the application
\begin{equation}
\label{AtHePe:32}
\begin{aligned}
\mathcal{A}_t:\H^4(\Er)\times\H^2(\Er)\times\mathrm{E}_2(t)&\rightarrow \mathrm{E}_0(t)\\
(u_0,u_1,g)&\mapsto u
\end{aligned}
\end{equation}
where \(u\) denotes the solution introduced in Proposition \ref{AtHePe:Prop1},
can be continuously extended to a function, still denoted by \(\mathcal{A}_t\):
\begin{equation*} 
\mathcal{A}_t: \H^2(\Er)\times\L^2(\Er)\times\L^2([0,t];\L^2(\Er))\rightarrow \mathrm{E}_1(t).
\end{equation*}
For \((u_0,u_1)\in\H^2(\Er)\times \L^2(\Er)\) fixed, we define
\begin{equation*}
\begin{aligned}
\mathcal{B}_t:\L^2([0,t];\L^2(\Er))&\rightarrow \mathrm{E}_1(t)\\
g&\mapsto \mathcal{A}_t(u_0,u_1,g).
\end{aligned}
\end{equation*}
If we denote by \(\mathcal{B}_t(g)\) by \(u\), then \(u\) still satisfies \eqref{AtHePe:3}. In the sequel,
we always assume that \((u_0,u_1)\in\H^2(\Er)\times\L^2(\Er)\).


\section{The approximated problem}
\label{Approx_Problem}


We mollify in system \eqref{AtHePe:2} the \(\delta\) function, mollify and truncate the 
nonlinear term.
Namely, we consider the following approximated problem:
\begin{subequations}
\label{AtHePe:4}
\begin{empheq}[left=\empheqlbrace]{alignat=1}
\label{AtHePe:4_a}
&u_{tt}+Fu_{xxxx}+Gu_t-((\varphi^R(u\star \theta')-{\nu p})(u\star\theta''))\star\theta+h)F=0,\\
\label{AtHePe:4_b}
&u(0,\cdot)=u_0\quad\text{and}\quad u_t(0,\cdot)=u_1,
\end{empheq}
\end{subequations}
where \(\star\) denotes the convolution with respect to \(x\), 
and \(h\in\C^0([0,T];\H^2(\Er))\) is a given function.
Here, for \(R>0\), \(\varphi^R\in\C^{\infty}(\Er)\) is an even, increasing function on \(\Er_+\)
such that \(\varphi^R(x)\eqldef x^2\) for \(\abs{x}\leq R\) and \(\varphi^{R}(x)\eqldef (R+1)^2\)
for \(\abs{x}\geq R+2\). Clearly, the equation \eqref{AtHePe:4_a} can formally be rewritten as follows:
\begin{equation*}
u_{tt}+Fu_{xxxx}+G u_t=\mathcal{C}_{t,R}(u),
\end{equation*}
where \(\mathcal{C}_{t,R}\) is defined by
\begin{equation*}
\begin{aligned}
\mathcal{C}_{t,R}: \mathrm{E}_1(t)&\rightarrow \mathrm{E}_2(t),\\
v&\mapsto \mathcal{C}_{t,R}(v)\eqldef \bigl([(\varphi^R(v\star \theta')-{\nu p})(v\star\theta'')]\star\theta+h\bigr)F
\end{aligned}
\end{equation*}
The fact that  \(\mathcal{C}_{t,R}\)  is well defined with values in
\(E_2(t)\) comes from usual convolution inequalities.
Similarly, we have
\begin{proposition}
\label{AtHePe:Prop2}
Let \(t\in[0,T]\). Then, 
for any \((v,\tilde v)\in\mathrm{E}_1(t)\times \mathrm{E}_1(t)\), the following inequality 
\begin{equation}
\label{AtHePe:5}
\begin{aligned}
&\norm[\L^2(0,t;\L^2(\Er))]{\mathcal{C}_{t,R}(v)-\mathcal{C}_{t,R}(\tilde v)}
\leq C_{T,R} \bigl(\norm[\L^2({[}0,t{]};\L^2(\Er))]{v-\tilde v}\\&+\norm[\L^2({[}0,t{]};\L^2(\Er))]{v_t-\tilde v_t}
+\norm[\C^0({[}0,t{]};\L^2(\Er))]{\tilde v}\norm[\L^2({[}0,t{]};\L^2(\Er))]{v-\tilde v}
\bigr)
\end{aligned}
\end{equation}
holds true.
\end{proposition}

\bproof
Let \((v,\tilde v)\in\mathrm{E}_1(t)\times \mathrm{E}_1(t)\). Hence, we have
\begin{equation*}
\begin{aligned}
&\norm[\L^2({[}0,t{]};\L^2(\Er))]{([(\varphi^R(v\star \theta')-{\nu p})(v\star\theta'')]\star\theta-
[(\varphi^R(\tilde v\star \theta')-{\nu p})(\tilde v\star\theta'')]\star\theta)F}
\\&\leq C_{T,R}\bigl((\norm[\L^2({[}0,t{]};\L^2(\Er))]{\varphi^R(v\star \theta')-{\nu p})(v-\tilde v)\star\theta''}
\\&+\norm[\L^2({[}0,t{]};\L^2(\Er))]{(\varphi^R(v\star\theta')-\varphi^R(\tilde v\star\theta'))(\tilde v\star\theta'')}
\bigr).
\end{aligned}
\end{equation*}
Since \(\varphi^{R}\) is bounded and globally Lipschitz and by convolution inequalities, it comes that
\begin{equation*}
\begin{aligned}
&\norm[\L^2({[}0,t{]};\L^2(\Er))]{([(\varphi^R(v\star \theta')-{\nu p})(v\star\theta'')]\star\theta-
[(\varphi^R(\tilde v\star \theta')-{\nu p})(\tilde v\star\theta'')]\star\theta)F}
\\&\leq 
C_{T,R}
\bigl(
\norm[\L^2({[}0,t{]};\L^2(\Er))]{v-\tilde v}
+\norm[\C^0({[}0,t{]};\L^2(\Er))]{\tilde v}\norm[\L^2({[}0,t{]};\L^2(\Er))]{v-\tilde v}
\bigr),
\end{aligned}
\end{equation*}
which proves the result.
\eproof

We are now looking for the fixed points of application:
\begin{equation*}
\mathcal{B}_T\mathcal{C}_{T,R}:\mathrm{E}_1(T)\rightarrow \mathrm{E}_1(T).
\end{equation*}

\begin{proposition}
\label{AtHePe:Prop3}
The application \(\mathcal{B}_T\mathcal{C}_{T,R}\) admits a fixed point.
\end{proposition}

\bproof
We use the Picard fixed point theorem. Due to Proposition \ref{AtHePe:Prop2} \(\emph{(i)}\), 
we have mainly to bound \(\norm[\C^0({[}0,t{]};\L^2(\Er))]{\tilde v}\).

\vspace{0.5em}
\noindent\textbf{(a) An invariant set.}
\vspace{0.5em}

\noindent First, let us assume that \(v\in\mathrm{E}_1(T)\) with \((v(0,\cdot),v_t(0,\cdot))=(u_0,u_1)\) 
and \(\tilde v=0\).  
It follows from \eqref{AtHePe:3} and \eqref{AtHePe:5} that
\begin{equation}
\label{AtHePe:6}
\begin{aligned}
&\norm[\H^2(\Er)]{(\mathcal{B}_T\mathcal{C}_{T,R}(v))(t)}
+{\|}{\partial_t}(\mathcal{B}_T\mathcal{C}_{T,R}(v))(t){\|}_{\L^2(\Er)}
\\&\leq 
C_{T,R}
\bigl(\norm[\H^2(\Er)]{u_0}+\norm[\L^2(\Er)]{u_1}+\norm[\L^2({[}0,T{]};\L^2(\Er))]{\mathcal{C}_{T,R}(v)}\bigr)
\\&\leq 
C_{T,R}
\bigl(\norm[\H^2(\Er)]{u_0}+\norm[\L^2(\Er)]{u_1}+\norm[\L^2({[}0,T{]};\L^2(\Er))]{\mathcal{C}_{T,R}(0)}
\\&\hspace{1em} +
C_{T,R}
\bigl(\norm[\L^2({[}0,t{]};\L^2(\Er))]{v}+\norm[\L^2({[}0,t{]};\L^2(\Er))]{v_t}\bigr)
\bigr).
\end{aligned}
\end{equation}
For \(f\in \C^0({[}0,t{]};\L^2(\Er))\), we define
\begin{equation*}
\norm[\lambda,t]{f}\eqldef \sup_{0\leq s\leq t}\bigl(\mathrm{e}^{-\lambda s}\norm[\L^2(\Er)]{f(s)}\bigr)
\end{equation*}
with \(\lambda>0\) and for \(g\in\mathrm{E}_1(t)\), we define
\begin{equation*}
\norm[1,\lambda,t]{g}\eqldef \sup_{0\leq s\leq t}\bigl(\mathrm{e}^{-\lambda s}\bigl(\norm[\H^2(\Er)]{g(s)}+\norm[\L^2(\Er)]{g_t(s)}\bigr)\bigr).
\end{equation*}
It comes from \eqref{AtHePe:6} that
\begin{equation}
\label{AtHePe:7}
\begin{aligned}
&\mathrm{e}^{-\lambda t}\bigl(\norm[\H^2(\Er)]{(\mathcal{B}_T\mathcal{C}_{T,R}(v))(t)}
+{\|}{\partial_t}(\mathcal{B}_T\mathcal{C}_{T,R}(v))(t){\|}_{\L^2(\Er)}\bigr)
\\&\leq 
C_{T,R}
\mathrm{e}^{-\lambda t}\bigl(\norm[\H^2(\Er)]{u_0}+\norm[\L^2(\Er)]{u_1}+\norm[\L^2({[}0,T{]};\L^2(\Er))]{\mathcal{C}_{T,R}(0)}\bigr)
\\&\hspace{1em}+
C_{T,R}
\Bigl(\int_0^t \mathrm{e}^{-2\lambda t}\mathrm{e}^{2\lambda s}\mathrm{e}^{-2\lambda s}\bigl(\norm[\L^2(\Er)]{v(s)}^2+\norm[\L^2(\Er)]{v_t(s)}^2\bigr)
\dd s\Bigr)^{1/2}
\\&\leq 
C_{T,R}
\mathrm{e}^{-\lambda t}\bigl(\norm[\H^2(\Er)]{u_0}+\norm[\L^2(\Er)]{u_1}+\norm[\L^2({[}0,T{]};\L^2(\Er))]{\mathcal{C}_{T,R}(0)}\bigr)
\\&\hspace{1em}+\frac{
C_{T,R}
}{\sqrt{\lambda}}\norm[1,\lambda,t]{v}.
\end{aligned}
\end{equation}
Note that \eqref{AtHePe:7} can be rewritten for any \(s\in [0,t]\)
in place of \(t\). Furthermore, in this new inequality, \(\norm[1,\lambda,s]{\,\cdot\,}\)
and \(\norm[\L^2(0,s;\L^2(\Er))]{\,\cdot\,}\) are respectively smaller that \(\norm[1,\lambda,t]{\,\cdot\,}\)
and \(\norm[\L^2(0,t;\L^2(\Er))]{\,\cdot\,}\). Hence, we get \eqref{AtHePe:7} with \(s\) in place of \(t\) in the left hand side. Taking \(t=T\), we obtain
\begin{equation}
\label{AtHePe:8}
\norm[1,\lambda,T]{\mathcal{B}_T\mathcal{C}_{T,R}(v)}
\leq M_R+\frac{C_{T,R}
}{\sqrt{\lambda}}\norm[1,\lambda,T]{v},
\end{equation}
where \(M_R\eqldef 
C_{T,R}
\bigl(\norm[\H^2(\Er)]{u_0}+\norm[\L^2(\Er)]{u_1}+\norm[\L^2({[}0,T{]};\L^2(\Er))]{\mathcal{C}_{T,R}(0)}\bigr)\).
Let \(\lambda=4
C_{T,R}^2\) and define \(\beta(0,2M_R)\eqldef\{u\in\mathrm{E}_1(T):\norm[1,\lambda,T]{u}\leq 2M_R\}\). Then, we may deduce from \eqref{AtHePe:8} that
\(\mathcal{B}_T\mathcal{C}_{T,R}(\beta(0,2M_R))\subset \beta(0,2M_R)\).

\vspace{0.5em}
\noindent\textbf{(b) Contraction of  \(\mathcal{B}_T\mathcal{C}_{T,R}\).} 
\vspace{0.5em}

Let \((v,\tilde v)\in \beta(0,2M_R)\times \beta(0,2M_R)\) with \(v(0,\cdot)=\tilde v(0,\cdot)=u_0\) and  \(v_t(0,\cdot)=\tilde v_t(0,\cdot)=u_1\).
Since \(v-\tilde v\) has \(0\) initial data, we have (see \eqref{AtHePe:3})
\begin{equation*}
\begin{aligned}
&\norm[\H^2(\Er)]{(\mathcal{B}_T\mathcal{C}_{T,R}(v))(t){-}(\mathcal{B}_T\mathcal{C}_{T,R}(\tilde v))(t)}
{+}{\|}\partial_t((\mathcal{B}_T\mathcal{C}_{T,R}(v))(t){-}(\mathcal{B}_T\mathcal{C}_{T,R}(\tilde v))(t))
{\|}_{\L^2(\Er)}
\\&
\leq 
C_{T,R}
\norm[\L^2({[}0,t{]};\L^2(\Er))]{\mathcal{C}_{T,R}(v)-\mathcal{C}_{T,R}(\tilde v)}
\\&\leq 
C_{T,R}
\bigl(\norm[\L^2({[}0,T{]};\L^2(\Er))]{v-\tilde v}+\norm[\L^2({[}0,T{]};\L^2(\Er))]{v_t-\tilde v_t}
\\&\hspace{1em}+\norm[\C^0({[}0,T{]};\L^2(\Er))]{\tilde v}\norm[\L^2({[}0,T{]};\L^2(\Er))]{v-\tilde v}
\bigr),
\end{aligned}
\end{equation*}
Since \(\tilde v\in\beta(0,2M_R)\), we have \(\norm[\C^0({[}0,t{]};\L^2(\Er))]{\tilde v}\leq 2M_R\mathrm{e}^{\lambda T}\). It follows that
\begin{equation*}
\begin{aligned}
&
\norm[\H^2(\Er)]{(\mathcal{B}_T\mathcal{C}_{T,R}(v))(t){-}(\mathcal{B}_T\mathcal{C}_{T,R}(\tilde v))(t)}
{+}{\|}\partial_t((\mathcal{B}_T\mathcal{C}_{T,R}(v))(t){-}(\mathcal{B}_T\mathcal{C}_{T,R}(\tilde v))(t))
{\|}_{\L^2(\Er)}
\\&
\leq 
C_{T,R}
\bigl(\norm[\L^2({[}0,T{]};\L^2(\Er))]{v-\tilde v}+\norm[\L^2({[}0,T{]};\L^2(\Er))]{v_t-\tilde v_t}\bigr).
\end{aligned}
\end{equation*}
Introducing a suitable \(\mu>0\), we derive as before the following inequality
\begin{equation*}
\norm[1,\mu,T]{\mathcal{B}_T\mathcal{C}_{T,R}(v)-\mathcal{B}_T\mathcal{C}_{T,R}(\tilde v)}
\leq \frac12 \norm[1,\mu,T]{v-\tilde v},
\end{equation*}
which proves the Proposition.
\eproof

An immediate consequence of Proposition \ref{AtHePe:Prop2} is the following Corollary:

\begin{corollary}
\label{AtHePe:Cor1}
Assume that \(u_0\in\H^4(\Er)\), \(u_1\in\H^2(\Er)\). Let \(u^R\in\mathrm{E}_1(T)\) such that \(\mathcal{B}_T\mathcal{C}_{T,R}(u^R)=u^R\).
Then, \(u^R\in\mathrm{E}_0(T)\). Moreover,
\begin{subequations}
\label{AtHePe:9}
\begin{empheq}[left=\empheqlbrace]{alignat=1}
\label{AtHePe:9_a}
&u_{tt}^R+Fu_{xxxx}^R+Gu_{t}^R-F((\varphi^R(u^R\star \theta')-{\nu p})(u^R\star\theta''))\star\theta+h)=0,\\
\label{AtHePe:9_b}
&u^R(0,\cdot)=u_0\quad\text{and}\quad u_t^{R}(0,\cdot)=u_1.
\end{empheq}
\end{subequations}
\end{corollary}

\bproof
We have \(\mathcal{C}_{T,R}(\mathrm{E}_1(T))\subset\mathrm{E}_2(T)\). 
Hence, by 
\eqref{AtHePe:32}, we get \(\mathcal{B}_T\mathcal{C}_{T,R}(u^R)\in\mathrm{E}_0(T)\), which is \(u^{R}\in\E_0(T)\).
\eproof

Let \(\epsilon>0\) and \(\theta^{\epsilon}(x)\eqldef \frac1{\epsilon}\theta\bigl(\frac{x}{\epsilon}\bigr)\) with \(x\in\Er\).
We assume now that \(u_0\in\H^2(\Er)\) and \(u_1\in\L^2(\Er)\cap\L^{\infty}(\Er)\) and set \(u^{\epsilon}_0\eqldef u_0\star\theta^{\epsilon}\),
\(v^{\epsilon}_0\eqldef u_1\star\theta^{\epsilon}\) and \(h^{\epsilon}\eqldef -(f+\theta^{\epsilon}(\cdot-\zeta(t))P(t))\). 
Clearly, we have \(u^{\epsilon}_0\in\H^{4}(\Er)\) and  \(u^{\epsilon}_1\in\H^{2}(\Er)\).
We apply the Corollary \ref{AtHePe:Cor1} with  \((u^{\epsilon}_0,u^{\epsilon}_1,\theta^{\epsilon},h^{\epsilon})\) in place of \((u_0,u_1,\theta,h)\).
This provides
a function \(u^{R,\epsilon}\) solution of the following problem:  
\begin{subequations}
\label{AtHePe:10}
\begin{empheq}[left=\empheqlbrace]{alignat=1}
\label{AtHePe:10_a}
&u^{R,\epsilon}_{tt}+\frac{\dd}{\dd t}(\theta^{\epsilon}(\cdot-\zeta(t))u^{R,\epsilon}_{t})
+u^{R,\epsilon}_{xxxx}\\&\notag
-((\varphi^R(u^{R,\epsilon}\star {\theta^{\epsilon}}')-{\nu p})(u^{R,\epsilon}\star{\theta^{\epsilon}}''))\star\theta^{\epsilon}+h^{\epsilon})=0,\\
\label{AtHePe:10_b}
&u^{R,\epsilon}(0,\cdot)=u^{\epsilon}_0\quad\text{and}\quad u^{R,\epsilon}_t(0,\cdot)=u^{\epsilon}_1.
\end{empheq}
\end{subequations}
In the sequel, we set  \(Q_t\eqldef [0,t]\times\Er\) where \(t\in[0,T]\). We have the following (uniform in \(R\) and \(\epsilon\)) estimates:

\begin{proposition}
\label{AtHePe:Prop4}
Assume that \(u_0\in\H^2(\Er)\), \(u_1\in\L^2(\Er)\cap\L^{\infty}(\Er)\). Then, there exists a constant \(C>0\), independent of 
\(R\) and \(\epsilon\), such that
\begin{equation*}
\sup_{0\leq t\leq T}\bigl(\norm[\H^2(\Er)]{u^{R,\epsilon}(t)}+\norm[\L^2(\Er)]{u^{R,\epsilon}_t(t)}\bigr)
\leq C.
\end{equation*}
\end{proposition}

\bproof
We define \(\frac{\partial}{\partial \tau}\eqldef \frac{\partial}{\partial t}+\dot\zeta(t)\frac{\partial}{\partial x}\). Then, we multiply \eqref{AtHePe:10_a}
by \(u^{R,\epsilon}_{\tau}\) and we integrate this result over \([0,t]\times\Er\) with \(t\in[0,T]\). Since \(u^{R,\epsilon}\in \mathrm{E}_0(T)\), the boundary 
terms in \(x=\pm\infty\) vanish. We observe that the term \(u^{R,\epsilon}_{tt}\) in \eqref{AtHePe:10} gives
\begin{equation*} 
\begin{aligned}
A_1(t)\eqldef\, & \frac12\norm[\L^2(\Er)]{u^{R,\epsilon}_{t}(t)}^2-\frac12\norm[\L^2(\Er)]{u_1}^2+\int_{Q_t}\dot\zeta u^{R,\epsilon}_{tt}u^{R,\epsilon}_{x} \dd s\dd x
\\ \geq\, &\frac12\norm[\L^2(\Er)]{u^{R,\epsilon}_{t}(t)}^2-\frac12\norm[\L^2(\Er)]{u_1}^2-\frac12\int_{Q_t}\partial_x(\dot \zeta(u^{R,\epsilon}_{t})^2)\dd s\dd x
\\&-\int_{Q_t}\ddot\zeta u^{R,\epsilon}_{x} u^{R,\epsilon}_t\dd s\dd x+\int_{\Er}\bigl{[}\dot\zeta u^{R,\epsilon}_t u^{R,\epsilon}_x\bigr{]}_0^t\dd x
\\\geq\, & \Bigl(\frac12-\eta\Bigr)\norm[\L^2(\Er)]{u^{R,\epsilon}_t(t)}^2-C\int_0^t\bigl
(\norm[\L^2(\Er)]{u^{R,\epsilon}_t(s)}^2+\norm[\L^2(\Er)]{u^{R,\epsilon}_x(s)}^2\bigr)\dd s
\\&-C_{\eta}\norm[\L^2(\Er)]{u^{R,\epsilon}_x(t)}^2-C\bigl(\norm[\H^2(\Er)]{u_0}^2+\norm[\L^2(\Er)]{u_1}^2\bigr).
\end{aligned}
\end{equation*}
Concerning the term \(\frac{\dd}{\dd t}(\theta_{\epsilon}(x-\zeta(t))u^{R,\epsilon}_t)\), we have
\begin{equation*}
\begin{aligned}
A_2(t)\eqldef& \int_{Q_t}
\frac{\dd}{\dd s}(\theta^{\epsilon}(x-\zeta(s))u^{R,\epsilon}_t) u^{R,\epsilon}_{\tau}\dd s\dd x
\\=&-\int_{Q_t}\theta^{\epsilon}(x-\zeta(s))u^{R,\epsilon}_t\partial_{\tau}\partial_{t}u^{R,\epsilon}+\ddot\zeta u^{R,\epsilon}_x)\dd s\dd x
+\int_{\Er}\bigl[\theta^{\epsilon}(x-\zeta(s))u_t^{R,\epsilon}u^{R,\epsilon}_{\tau}\bigr]_0^t\dd x
\\=&-\frac12\int_{Q_t} \partial_{\tau}\bigl(\theta^{\epsilon}(x-\zeta(s))(u^{R,\epsilon}_t)^2\bigr)\dd s\dd x
-\int_{Q_t}\theta^{\epsilon}(x-\zeta(t))\ddot \zeta u^{R,\epsilon}_tu^{R,\epsilon}_x\dd s\dd x
\\&+\int_{\Er}\theta^{\epsilon}(x-\zeta(t))\bigl((u^{R,\epsilon}_t)^2+\dot\zeta(t) u^{R,\epsilon}_t(t,x)u^{R,\epsilon}_x(t,x)\bigr)\dd x
\\&-\int_{\Er}\theta^{\epsilon}(x-\zeta(0))\bigl((u^{\epsilon}_1)^2+\dot\zeta(0)u^{\epsilon}_1 {u^{\epsilon}_{0}}'\bigr)\dd x
\\\geq& -\frac12\int_{\Er}\bigl(\theta^{\epsilon}(x-\zeta(t))(u^{R,\epsilon}_t(t,x))^2-\theta^{\epsilon}(x-\zeta(0))(u^{\epsilon}_1)^2\bigr)\dd x
\\&-\int_{Q_t}\theta^{\epsilon}(x-\zeta(s))(u^{R,\epsilon}_t)^2\dd t\dd x
\\&-C\int_{Q_t}\theta^{\epsilon}(x-\zeta(s))\bigl(\norm[\L^2(\Er)]{u^{R,\epsilon}(s)}^2
+\norm[\L^2(\Er)]{u^{R,\epsilon}_{xx}(s)}^2\bigr)\dd s\dd x
+\int_{\Er}\theta^{\epsilon}(x-\zeta(t))(u^{R,\epsilon}_t(t,x))^2\dd x
\\&-C_{\eta}\int_{\Er}\theta^{\epsilon}(x-\zeta(t))(u^{R,\epsilon}_x(t,x))^2\dd x-\eta\int_{\Er}\theta^{\epsilon}(x-\zeta(t))(u^{R,\epsilon}_t(t,x))^2\dd x
\\&-\int_{\Er}\theta^{\epsilon}(x-\zeta(0))(u^{\epsilon}_1)^2\dd x
-\int_{\Er}\theta_{\epsilon}(x-\zeta(0))\dot\zeta(0)u^{\epsilon}_1{u_0^{\epsilon}}'\dd x
\\\geq&\, \Bigl(\frac12-\eta\Bigr)\int_{\Er}\theta^{\epsilon}(x-\zeta(t))(u^{R,\epsilon}_t(t,x))^2\dd x 
-\int_{Q_t}\theta^{\epsilon}(x-\zeta(s))(u^{R,\epsilon}_t)^2\dd s\dd x
\\&-C\int_0^t\norm[\H^2(\Er)]{u^{R,\epsilon}(s)}^2\dd s -C_{\eta}\norm[\L^2(\Er)]{u^{R,\epsilon}(t)}^2
\\&-\eta\norm[\H^2(\Er)]{u^{R,\epsilon}(t)}^2-C\bigl(\norm[\L^{\infty}(\Er)]{u_1}^2+\norm[\H^2(\Er)]{u_0}^2
\bigr)
\end{aligned}
\end{equation*}
The term \(u^{R,\epsilon}_x\) provides
\begin{equation*}
\begin{aligned}
A_3(t)\eqldef& \int_{Q_t}u^{R,\epsilon}_{xxxx}u^{R,\epsilon}_{\tau}\dd s\dd x
\\\geq&\, \frac12\norm[\L^2(\Er)]{u^{R,\epsilon}_{xx}(t)}^2-\frac12\norm[\L^2]{u_0''}^2
+\frac12\int_0^t\dot\zeta(s)\int_{\Er}\partial_{x}(u^{R,\epsilon}_{xx})^2\dd s\dd x
\\=&\,\frac12\bigl(\norm[\L^2(\Er)]{u^{R,\epsilon}_{xx}(t)}^2-\norm[\H^2(\Er)]{u_0}^2\bigr).
\end{aligned}
\end{equation*}
In the sequel, we define \(\psi^R\) and \(\mu^R\) by \({\psi^R}'\eqldef \varphi^R\) and \({\mu^R}'\eqldef \psi^R\) 
such that \(\mu^R(0)=\psi^R(0)=0\). Since \(\varphi^R\) is even, so is \(\mu^R\). Furthermore, notice that \({\mu^R}''=\varphi^R\geq 0\).
Consequently, \(\mu^R\) is convex, even with \(\mu^R(0)=0\). It follows that \(\mu^R\geq 0\).
Now, the term 
\(-(\varphi^R(u^{R,\epsilon}_x\star \theta^{\epsilon})-{\nu p})(u^{R,\epsilon}_{xx}\star\theta^{\epsilon}))\star\theta^{\epsilon}\) provides
\begin{equation*}
\begin{aligned}
A_4(t)\eqldef& 
-\int_{Q_t}(\varphi^R(u^{R,\epsilon}_x\star \theta^{\epsilon})-{\nu p})(u^{R,\epsilon}_{xx}\star\theta_{\epsilon}))\partial_{\tau}(u^{R,\epsilon}\star\theta^{\epsilon})\dd s\dd x
\\=&\int_{Q_t}(\psi_R(u^{R,\epsilon}_x\star\theta^{\epsilon})\partial_{\tau}(u^{R,\epsilon}_x\star\theta^{\epsilon})\dd s\dd x+
\int_0^t\bigl[
\psi^R(u^{R,\epsilon}_x\star\theta^{\epsilon})\partial_{\tau}(u^{R,\epsilon}\star\theta^{\epsilon})
\bigr]_{-\infty}^{+\infty}\dd s
\\&+\int_{Q_t}{\nu p}(u^{R,\epsilon}_{xx}\star\theta^{\epsilon})\frac{\partial}{\partial \tau}(u^{R,\epsilon}\star\theta^{\epsilon})\dd s\dd x.
\end{aligned}
\end{equation*}
Since \(\abs{\psi^R(u^{R,\epsilon}\star{\theta^{\epsilon}}')}\leq C\abs{u^{R,\epsilon}\star{\theta^{\epsilon}}'}\leq C_{R,\epsilon}\)
and \(\lim_{\abs{x}\rightarrow +\infty}\bigl{|}\frac{\partial}{\partial\tau}(u^{R,\epsilon}\star\theta^{\epsilon})\bigr{|}=0\), we may deduce that
\(\int_0^t\bigl[
\psi^R(u^{R,\epsilon}_x\star\theta^{\epsilon})\partial_{\tau}(u^{R,\epsilon}_x\star\theta^{\epsilon})
\bigr]_{-\infty}^{+\infty}\dd t\) vanishes. Hence, we have
\begin{equation*}
\begin{aligned}
A_4(t)\geq& \int_{Q_t}\frac{\partial}{\partial\tau}(\mu^R(u^{R,\epsilon}_x\star\theta^{\epsilon}))\dd s\dd x-
C\int_0^t\bigl(\norm[\H^2(\Er)]{u^{R,\epsilon}(s)}^2+\norm[\L^2(\Er)]{u^{R,\epsilon}_t(s)}^2\bigr)\dd s
\\=&\int_{\Er}\bigl(\mu^R(u^{R,\epsilon}_x(t)\star\theta^{\epsilon}(t))-\mu^R(u^{R,\epsilon}_x(0,x)\star\theta^{\epsilon}(0))\bigr)\dd x
\\&-C\int_0^t\bigl(\norm[\H^2(\Er)]{u^{R,\epsilon}(s)}^2+\norm[\L^2(\Er)]{u^{R,\epsilon}_t(s)}^2\bigr)\dd s.
\end{aligned}
\end{equation*}
Since \(\mu^R(u^{R,\epsilon}_x(t)\star\theta_{\epsilon}(t))\geq 0\), it comes that
\begin{equation*}
\begin{aligned}
A_4(t)\geq&  -C\int_{\Er}\abs{{u^{\epsilon}_{0}}'}^4\dd x-C\int_0^t\bigl(\norm[\H^2(\Er)]{u^{R,\epsilon}(s)}^2+\norm[\L^2(\Er)]{u_t^{R,\epsilon}(s)}^2\bigr)\dd s
\\\geq& -C\norm[\H^2(\Er)]{u_0}^4-C\int_0^t\bigl(\norm[\H^2(\Er)]{u^{R,\epsilon}(s)}^2+\norm[\L^2(\Er)]{u^{R,\epsilon}_t(s)}^2\bigr)\dd s.
\end{aligned}
\end{equation*}
The term \(h_{\epsilon}\) provides
\begin{equation*}
\begin{aligned}
A_5(t)\eqldef& -\int_{Q_t} f u^{R,\epsilon}_\tau\dd s\dd x
-\int_{Q_t} P(s)\theta_{\epsilon}(x-\zeta(t))u^{R,\epsilon}_{\tau}\dd s\dd x\\
\geq& -C\int_{Q_t}f^2\dd s\dd x-C\int_{Q_t}\abs{u^{R,\epsilon}_t}^2\dd s\dd x-C\int_0^t\norm[\H^2(\Er)]{u^{R,\epsilon}(s)}^2\dd s\\&
+\int_{Q_t}P(s)\partial_{\tau}(\theta^{\epsilon}(x-\zeta(s)))u^{R,\epsilon}(s,x)\dd s\dd x-
\int_{\Er}\bigl[P(s)\theta^{\epsilon}(x-\zeta(s))u^{R,\epsilon}(s,x)\bigr]_0^t\dd x\\
\geq& -C\int_{Q_t} f^2-C\int_{Q_t}\abs{u^{R,\epsilon}_t}^2\dd s\dd x-C\int_0^t\norm[\H^2(\Er)]{u^{R,\epsilon}(s)}^2\dd s-C\norm[\H^2(\Er)]{u_0}^2
\\&-C_{\eta}\norm[\L^2(\Er)]{u^{R,\epsilon}(t)}^2-\eta \norm[\H^2(\Er)]{u^{R,\epsilon}(t)}^2.
\end{aligned}
\end{equation*}
Since \(\sum_{k=1}^5A_k=0\), the previous estimates lead to
\begin{equation}
\begin{aligned}
\label{AtHePe:11}
&\norm[\L^2(\Er)]{u^{R,\epsilon}_t(t)}^2+\norm[\L^2(\Er)]{u^{R,\epsilon}_{xx}(t)}^2
+\int_{\Er}\theta^{\epsilon}(x-\zeta(s))(u^{R,\epsilon}_t(t,x))^2\dd x
\\&\leq C_T\Bigl(\norm[\H^2(\Er)]{u_0}^2+\norm[\H^2(\Er)]{u_0}^4+\norm[\L^2(\Er)\cap\L^{\infty}(\Er)]{u_1}^2+
\int_0^t\norm[\H^2(\Er)]{u^{R,\epsilon}(s)}^2\dd s
\\&+\int_0^t\norm[\L^2(\Er)]{u^{R,\epsilon}_t(s)}^2\dd s
+\int_{Q_t}\theta^{\epsilon}(x-\zeta(s))(u^{R,\epsilon}_t(s,x))^2\dd s\dd x+1
\Bigr)\\&+2\eta\norm[\H^2(\Er)]{u^{R,\epsilon}(t)}^2+C_{\eta}\norm[\L^2(\Er)]{u^{R,\epsilon}(t)}^2.
\end{aligned}
\end{equation}
In order to handle the terms \(\eta\norm[\H^2(\Er)]{u^{R,\epsilon}(t)}^2\) and \(C_{\eta}\norm[\L^2(\Er)]{u^{R,\epsilon}(t)}^2\), we write
\begin{equation*}
u^{R,\epsilon}(t,x)=u_0(x)+\int_0^tu^{R,\epsilon}_t(s,x)\dd s,
\end{equation*}
which implies that
\begin{equation}
\label{AtHePe:12}
\norm[\L^2(\Er)]{u^{R,\epsilon}(t)}\leq C\Bigl(
\norm[\L^2(\Er)]{u_0}+
\int_0^t\norm[\L^2(\Er)]{u^{R,\epsilon}_t(s)}\dd s
\Bigr),
\end{equation}
and 
\begin{equation}
\label{AtHePe:12hh}
\norm[\L^2(\Er)]{u^{R,\epsilon}(t)}^2\leq C_T\Bigl(
\norm[\L^2(\Er)]{u_0}^2+
\int_0^t\norm[\L^2(\Er)]{u^{R,\epsilon}_t(s)}^2\dd s
\Bigr).
\end{equation}
It follows from \eqref{AtHePe:12hh} that
\begin{equation}
\label{AtHePe:12pp}
\norm[\H^2(\Er)]{u(t)}^2\leq C_T
\Bigl(\norm[\L^2(\Er)]{u_0}^2+
\int_0^t\norm[\L^2(\Er)]{u^{R,\epsilon}_t(s)}^2\dd s
+\norm[\L^2(\Er)]{u^{R,\epsilon}_{xx}(t)}^2
\Bigr).
\end{equation}
As a consequence of \eqref{AtHePe:12hh} and \eqref{AtHePe:12pp}, \eqref{AtHePe:11} remains valid with \(\norm[\H^2(\Er)]{u^{R,\epsilon}(t)}^2\)
in place of \(\norm[\L^2(\Er)]{u^{R,\epsilon}_{xx}(t)}^2\) in the left hand side, and \(\eta=0\) and \(C_{\eta}=0\)
in the right hand side. Finally, the Gr\"onwall's lemma gives the desired estimate.
\eproof

Now, from Proposition \ref{AtHePe:Prop4} and the Sobolev embeddings,  
there exists a constant \(C>0\), independent of \(R\) and \(\epsilon\), such that
\begin{equation*}
\norm[\L^{\infty}({[}0,T{]}\times\Er)]{u^{R,\epsilon}}\leq C.
\end{equation*}
Hence, \(\norm[\L^{\infty}({[}0,T{]}\times\Er)]{u^{R,\epsilon}\star{\theta^{\epsilon}}'}\leq \frac{C}\epsilon\).
It follows that for \(R=\frac{C}{\epsilon}\), we have \(\varphi_R(u^{R,\epsilon}\star{\theta^{\epsilon}}')=(u^{R,\epsilon}\star{\theta^{\epsilon}}')^2\).
Consequently, we replace \(\varphi^R(u^{R,\epsilon}\star{\theta^{\epsilon}}')\) by \((u^{R,\epsilon}\star{\theta^{\epsilon}}')^2\)  
in the identity \eqref{AtHePe:10_a} and drop the index \(R\) in \eqref{AtHePe:10_a}.  From \eqref{AtHePe:10}, we induce the following weak formulation:
\begin{equation}
\label{AtHePe:10faible}
\begin{cases}
\displaystyle{\int_0^T\int_{\Er}u_t^{\epsilon} v_t\dd x\dd t+\int_0^T\int_{\Omega} u_{xx}^{\epsilon}v_{xx}\dd x\dd t}
\\\noalign{\vspace*{\jot}} 
\displaystyle{+\frac13\int_0^T\int_{\Er}(
(u_x^{\epsilon}{\star}\theta^{\epsilon})^3(v_x{\star}\theta^{\epsilon}){+}\nu p(u_x^{\epsilon}{\star}\theta^{\epsilon})
(v_x{\star}\theta^{\epsilon}){+} (\nu p)_{x}(u_x^{\epsilon}{\star}\theta^{\epsilon}))\dd x\dd t}
\\\noalign{\vspace*{\jot}} 
\displaystyle{+\int_0^T\int_{\Er}\theta^{\epsilon}(x-\xi(t))u_t^{\epsilon} v_t\dd x\dd t
+\int_0^T\int_{\Er} fv\dd x\dd t+\int_0^TP(t)v(t,\xi(t))\dd t}
\\\noalign{\vspace*{\jot}} 
\displaystyle{+\int_{\Er} u_1^{\epsilon}(x)v(0,x)\dd x
+\int_{\Er}\theta^{\epsilon}(x-\xi(0))u_1^{\epsilon}(x) v(0,x)}\dd x=0,
\end{cases}
\end{equation}
for any \(v\in\C^2([0,T];\H^2(\Omega))\) with compact support in \([0,T)\times\Er\).


\section{Proof of Theorem \ref{AtHePe:thm1}}
\label{Passage_limit}


We derive from Proposition \ref{AtHePe:Prop4} the following convergences.

\begin{corollary}
\label{AtHePe:Cor2}
Assume that \(u_0\in\H^2(\Er)\), \(u_1\in\L^2(\Er)\cap\L^{\infty}(\Er)\). Then, there exists
a subsequence \(\{{u_{\epsilon_k}}\}_{k\in\En}\) of \(u_{\epsilon}\) such that 
\begin{enumerate}[(i)]

\item \(u^{\epsilon_k}, u^{\epsilon_k}_t, u^{\epsilon_k}_{xx}\underset{k\to +\infty}{\rightharpoonup} u, u_t, u_{xx}\text{ weakly in }\L^2(0,T;\L^2(\Er))\).

\item For any \(\alpha\in\,]0,2]\) and any compact set \(K\subset\Er\), we have
\begin{subequations}
\label{AtHePe:conv1-2}
\begin{align}
\label{AtHePe:conv1}
u^{\epsilon_k}&\underset{k\to +\infty}{\rightarrow} u\text{ strongly in }\C^0([0,T];\H^{2-\alpha}(K)),\\
\label{AtHePe:conv2}
u^{\epsilon_k}_x&\underset{k\to +\infty}{\rightarrow} u_x\text{ strongly in }\L^{\infty}([0,T]\times K).
\end{align}
\end{subequations}

\item For any 
\((a,b)\in\Er^2\) such that \(a<b\), \(T>0\), we have
\begin{subequations}
\label{AtHePe:13}
\begin{align}
\label{AtHePe:13a}
u^{\epsilon_k}_x\star\theta^{\epsilon_k}&\underset{k\to +\infty}{\rightarrow} u_x\text{ strongly in }\L^2([0,T]\times [a,b]),\\
\label{AtHePe:13b}
(u^{\epsilon_k}_x\star\theta^{\epsilon_k})^3&\underset{k\to +\infty}{\rightarrow} u_x^3\text{ strongly in }\L^2([0,T]\times [a,b]).
\end{align}
\end{subequations}

\end{enumerate}
\end{corollary}

\bproof
\(\emph{(i)}\) follows from Proposition \ref{AtHePe:Prop4}. Note that \eqref{AtHePe:conv1} comes from the Sobolev embeddings and
Aubin-Lions lemma while  \eqref{AtHePe:conv2}  follows from \eqref{AtHePe:conv1}  and 
\(\H^{1-\alpha}(K)\hookrightarrow\L^{\infty}(K)\) for \(\alpha\in\,\bigl{]}0,\frac12\bigr{[}\,\). It remains to prove \(\emph{(iii)}\).
Let \(t\in[0,T]\). We have
\begin{equation*}
\begin{aligned}
\norm[\L^2({[}a,b{]})]{u^{\epsilon_k}_x(t)\star\theta^{\epsilon_k}-u_x(t)}
&\leq 
\norm[\L^2({[}a,b{]})]{u^{\epsilon_k}_x(t)\star\theta^{\epsilon_k}-u_x(t)\star\theta^{\epsilon_k}}
+\norm[\L^2({[}a,b{]})]{u_x(t)\star\theta^{\epsilon_k}-u_x(t)}
\\&\leq 
\norm[\L^2({[}a,b{]})]{u^{\epsilon_k}_x(t)-u_x(t)}
+\norm[\L^2({[}a,b{]})]{u_x(t)\star\theta^{\epsilon_k}-u_x(t)}.
\end{aligned}
\end{equation*}
Since 
\begin{equation}
\label{AtHePe:14}
u^{\epsilon_k}_x\underset{k\to +\infty}{\rightarrow} u_x\text{ strongly in } \L^{\infty}([0,T]\times[a,b]),
\end{equation}
we have 
\(\norm[\L^2({[}a,b{]})]{u^{\epsilon_k}_x(t)-u_x(t)}\underset{k\to +\infty}{\rightarrow} 0\).
On the other hand, \(\norm[\L^2({[}a,b{]})]{u_x(t)\star\theta^{\epsilon_k}-u_x(t)}\underset{k\to +\infty}{\rightarrow} 0\).
Hence, we get  
\(\norm[\L^2({[}a,b{]})]{u^{\epsilon_k}_x(t)\star\theta^{\epsilon_k}-u_x(t)}
\underset{k\to +\infty}
{\rightarrow} 0\).
Besides, observe that
\begin{equation*}
\norm[\L^2({[}a,b{]})]{u^{\epsilon_k}_x(t)\star\theta^{\epsilon_k}-u_x(t)}
\leq 
\norm[\L^2({[}a,b{]})]{u^{\epsilon_k}_x(t)}+
\norm[\L^2({[}a,b{]})]{u_x(t)}
{\stackrel{\eqref{AtHePe:14}}{\leq}} C.
\end{equation*}
Finally, we deduce from the dominated convergence theorem that \eqref{AtHePe:13a} holds. It remains to prove \eqref{AtHePe:13b}.
To this aim, we may assume that \(\{u^{\epsilon_k}_x\}_{k\in\En}\) is bounded on 
\(\L^2([0,T]\times [a-1,b+1])\).
Hence, we have
\begin{equation*}
\begin{aligned}
&\norm[\L^2({[}a,b{]})]{(u^{\epsilon_k}_x\star\theta^{\epsilon_k})^3-u_x^3}
\\&\leq \norm[\L^2({[}a,b{]})]{(u^{\epsilon_k}_x\star\theta^{\epsilon_k}-u_x)((u^{\epsilon_k}_x\star\theta^{\epsilon_k})^2+(u^{\epsilon_k}_x\star\theta^{\epsilon_k})u_x+u_x^2)}
\\&{\stackrel{\eqref{AtHePe:14}}{\leq}} 
C\norm[\L^2({[}a,b{]})]{u^{\epsilon_k}_x\star\theta^{\epsilon_k}-u_x}.
\end{aligned}
\end{equation*}
The result follows from \eqref{AtHePe:13a}. 
\eproof

\begin{proposition}
\label{AtHePe:Prop5}
Assume that \(u_0\in\H^2(\Er)\) and \(u_1\in\H^1(\Er)\). 
Let \(v\in\C^2([0,T];\H^2(\Er))\) with a compact support \([0,T[\,\times\Er\).
Then, we have
\begin{subequations}
\label{AtHePe:Prop5e}
\begin{align}
\label{AtHePe:Prop5_1}
&\int_{Q_T}u_t^{\epsilon_k}v_t\dd t\dd x\underset{k\to +\infty}{\rightarrow} \int_{Q_t}u_tv_t\dd t\dd x,\\
\label{AtHePe:Prop5_2}
&\int_{Q_T}u_{xx}^{\epsilon_k}v_{xx}\dd t\dd x\underset{k\to +\infty}{\rightarrow} \int_{Q_t}u_{xx}v_{xx}\dd t\dd x,\\
\label{AtHePe:Prop5_3}
&\int_{Q_T}\Bigl(-\frac13(u^{\epsilon_k}_x\star\theta^{\epsilon_k})^3(v_x\star\theta^{\epsilon_k})
+{\nu p}(u^{\epsilon_k}_x\star\theta^{\epsilon_k})(v_x\star\theta^{\epsilon_k})
\\\notag&+(\nu p)_x(u^{\epsilon_k}_x\star\theta^{\epsilon_k})(v\star\theta^{\epsilon_k})
\Bigr)\dd t\dd x
\underset{k\to +\infty}{\rightarrow}
\int_{Q_T}\Bigl(-\frac13 u_x^3v_x-{\nu p} u_{xx} v\Bigr)\dd t\dd x,\\&
\label{AtHePe:Prop5_4}
\int_{Q_T}(f+\theta^{\epsilon}(x-\zeta(t))P(t))v\dd t\dd x
\underset{\epsilon\to 0}{\rightarrow}
\int_{Q_T} f v\dd x\dd t+\int_0^TP(t)v(t,\zeta(t))\dd t,\\&
\label{AtHePe:Prop5_5}
\int_{\Er}u_1^{\epsilon_k}(x)v(0,x)\dd x+\int_{\Er}\theta^{\epsilon_k}(x-\zeta(0))u_1^{\epsilon_k}(x)v(0,x)\dd x
\\&\notag\underset{k\to +\infty}{\rightarrow}
\int_{\Er} u_1(x)v(0,x)\dd x+u_1(\zeta(0))v(0,\zeta(0)),\\&
\label{AtHePe:Prop5_6}
\int_{Q_T}\theta^{\epsilon_k}(x-\zeta(t)) u^{\epsilon_k}_tv_t\dd t\dd x
\underset{k\to +\infty}{\rightarrow}
-u_0(\zeta(0))v_t(0,\zeta(0))
\\&\notag-\int_0^Tu(t,\zeta(t))v_{tt}(t,\zeta(t))\dd t\dd x
-\int_0^T\dot\zeta(t) u_x(t,\zeta(t))v_t(t,\zeta(t))\dd t
\\\notag&-\int_0^T\dot\zeta(t)u(t,\zeta(t)v_{xt}(t,\zeta(t))\dd t.
\end{align}
\end{subequations}
\end{proposition}

\bproof
We observe that \eqref{AtHePe:Prop5_1} and  \eqref{AtHePe:Prop5_2} follow from Corollary \ref{AtHePe:Cor2}\(\emph{(i)}\).
Notice that all the functions \(v_{xx}\), \(v\star\theta^{\epsilon_k}\) and \(v_x\star\theta^{\epsilon_k}\)
have their supports included in a fixed bounded product \([0,T]\times[a,b]\). Hence, we may deduce from \eqref{AtHePe:13}
that
\begin{equation*}
\begin{aligned}
&\int_{Q_T}\Bigl(-\frac13(u^{\epsilon_k}_x{\star}\theta^{\epsilon_k})^3(v_x{\star}\theta^{\epsilon_k})
+{\nu p}(u^{\epsilon_k}_x{\star}\theta^{\epsilon_k})(v_x{\star}\theta^{\epsilon_k})
+(\nu p)_x(u^{\epsilon_k}_x{\star}\theta^{\epsilon_k})(v{\star}\theta^{\epsilon_k})
\Bigr)\dd t\dd x
\\&\underset{k\to +\infty}{\rightarrow}
\int_{Q_T}\Bigl(-\frac13 u_x^3v_x+{\nu p}u_xv_x+(\nu p v)_xu_x\Bigr)\dd t\dd x
=\int_{Q_T}\Bigl(-\frac13 u_x^3v_x-{\nu p} u_{xx} v\Bigr)\dd t\dd x,
\end{aligned}
\end{equation*}
which proves \eqref{AtHePe:Prop5_3}.
Property \eqref{AtHePe:Prop5_4} is straightforward. Next, \(u_1\in\H^1(\Er)\subset(\L^2(\Er)\cap\L^{\infty}(\Er))\)
and \(u_1^{\epsilon}=u_1\star\theta^{\epsilon}\). Hence, we have
\begin{equation*}
u_1^{\epsilon}\underset{\epsilon\to 0}{\rightarrow}
u_1\text{ strongly in }\L^2(\Er)\cap\L^{\infty}(\Er),
\end{equation*}
and \eqref{AtHePe:Prop5_5} follows. In order to prove \eqref{AtHePe:Prop5_6}, we first establish the following result:
\begin{lemma}
\label{AtHePe:Lem2}
Let \(u_0\in\H^2(\Er)\) and \(u_1\in\H^1(\Er)\). Then, for any 
\(v\in\C^2([0,T];\H^2(\Er))\) with a compact support in \([0,T[\,\times \Er\), we have
\begin{equation*}
\begin{aligned}
&\int_{Q_T}\theta^{\epsilon_k}(x-\zeta(t))u^{\epsilon_k}_tv_t\dd t\dd x
=-\int_{\Er}\theta^{\epsilon_k}(x-\zeta(0))u^{\epsilon_k}_0 v_t(0,x)\dd x\\&
-\int_{Q_T}\theta^{\epsilon_k}(x-\zeta(t))u^{\epsilon_k} v_{tt}\dd t\dd x
-\int_{Q_T}\dot\zeta(t)\theta^{\epsilon_k}(x-\zeta(t)) u^{\epsilon_k}_xv_t\dd x\dd t
\\&-\int_{Q_T}\dot\zeta(t)\theta^{\epsilon_k}(x-\zeta(t)) u^{\epsilon_k}v_{xt}\dd x\dd t.
\end{aligned}
\end{equation*}
\end{lemma}

\bproof
Observe that
\begin{equation*}
\begin{aligned}
&\int_{Q_T}\theta^{\epsilon_k}(x-\zeta(t))u^{\epsilon_k}_tv_t\dd t\dd x
\\=&-\int_{Q_T}(-\dot\zeta(t){\theta^{\epsilon_k}}'(x-\zeta(t))v_t+\theta^{\epsilon_k}(x-\zeta(t))v_{tt})
u^{\epsilon_k}\dd t\dd x
\\&+\int_{\Er}\bigl[\theta^{\epsilon_k}(x-\zeta(t))u^{\epsilon_k} v_t\bigr]_0^T\dd x
\\=&
-\int_{\Er}\theta^{\epsilon_k}(x-\zeta(0))u^{\epsilon_k}_0v_t(0,x)\dd x
-\int_{Q_T}\theta^{\epsilon_k}(x-\zeta(t))u^{\epsilon_k}v_{tt}\dd t\dd x
\\&-\int_{Q_T}\dot\zeta(t)\theta^{\epsilon_k}(x-\zeta(t))\partial_x(u^{\epsilon_k}v_t)\dd t\dd x,
\end{aligned}
\end{equation*}
and the desired result follows.
\eproof

\noindent We now pass to the limit in the various terms appearing in the equality of Lemma \ref{AtHePe:Lem2}. 

\vspace{0.5em}
\noindent {\emph{End of the proof of \eqref{AtHePe:Prop5_6}}}.
We only treat the case of 
\(\int_{Q_T}\dot\zeta(t)\theta^{\epsilon_k}(x-\zeta(t))u^{\epsilon_k}_xv_t\dd t\dd x\),
which we replace by \(\int_0^T\int_a^b\dot\zeta(t)\theta^{\epsilon_k}(x-\zeta(t))u^{\epsilon_k}_xv_t\dd t\dd x\).
The other terms in Lemma \ref{AtHePe:Lem2} can be handled similarly. Notice that for \(\alpha\,\in\,\bigl{]}0,\frac12\bigr{[}\,\),
\(\H^{1-\alpha}([a,b])\) is a multiplicative algebra. Since \(u_{\epsilon_k,x}\in\C^0([0,T];\H^{1-\alpha}([a,b]))\)
and \(v_t\in\C^1([0,T];\H^1([a,b]))\hookrightarrow \C^0([0,T];\H^{1-\alpha}([a,b]))\), it follows that
\(u^{\epsilon_k}_x v_t\in  \C^0([0,T];\H^{1-\alpha}([a,b]))\) for any \(\alpha\in ,\bigl{]}0,\frac12\bigr{[}\,\).
Notice also, since \(\{u^{\epsilon_k}_x\}_{n\in\En^{\star}}\) is bounded in \(\L^{\infty}([0,T]\times[a,b])\) 
and \(v_t\in \C^1([0,T];\H^1([a,b]))\hookrightarrow \L^{\infty}([0,T]\times[a,b])\), that
\(\{u^{\epsilon_k}_xv_t\}_{n\in\En^{\star}}\) is bounded in \( \L^{\infty}([0,T]\times[a,b])\).
Now, we write 
\begin{equation}
\label{AtHePe:15}
\begin{aligned}
&\int_{Q_T}\dot\zeta(t)\theta^{\epsilon_k}(x-\zeta(t))u^{\epsilon_k}_xv_t\dd t\dd x
\\&=\int_{Q_T}\dot\zeta (t)\theta^{\epsilon_k}(x-\zeta(t))(u^{\epsilon_k}_x-u_x)v_t\dd t\dd x
+\int_{Q_T}\dot\zeta(t)\theta^{\epsilon_k}(x-\zeta(t))u_xv_t\dd t\dd x.
\end{aligned}
\end{equation}
Since \(u_xv_t\in \C^0([0,T];\H^{1-\alpha}([a,b]))\) and
\(\H^{1-\alpha}([a,b])\hookrightarrow\C^0([a,b])\) with \(\alpha\in\bigl{]}0,\frac12\bigr{[}\,\),
we have 
\begin{equation*}
\int_{\Er}\theta^{\epsilon_k}(x-\zeta(t))u_x(t,x)v_t(t,x)\dd x
\underset{k\to +\infty}{\rightarrow}
u_x(t,\zeta(t))v_t(t,\zeta(t))
\end{equation*}
for any \(t\in[0,T]\) and since 
\begin{equation*}
\Bigl{|}
\dot\zeta(t)\int_a^b\theta^{\epsilon_k}(x-\zeta(t))u_x(t,x)v_t(t,x)\dd x
\Bigr{|}\leq \norm[\L^{\infty}(0,T)]{\dot\zeta}\norm[\L^{\infty}({[}0,T{]}\times {[}a,b{]})]{u_x v_t},
\end{equation*}
it follows from the dominated convergence theorem that
\begin{equation*}
\int_{Q_T}\dot\zeta\theta^{\epsilon_k}(x-\zeta(t)) u_xv_t\dd t\dd x
\underset{k\to +\infty}{\rightarrow}
\int_0^T\dot\zeta(t) u_x(t,\zeta(t))v_t(t,\zeta(t)).
\end{equation*}
For the other term in \eqref{AtHePe:15}, we have
\begin{equation*}
\begin{aligned}
\Bigl{|}
&\int_{Q_T}
\dot\zeta(t)\theta^{\epsilon_k}(x-\zeta(t))(u_x^{\epsilon_k}(t,x)-u_x(t,x))v_t(t,x)\dd t\dd x
\Bigr{|}
\\&\leq
C_T \norm[\L^{\infty}(0,T)]{\dot\zeta}\norm[\L^{\infty}({[}0,T{]}\times {[}a,b{]})]{u^{\epsilon_k}_x-u_x}
\underset{k\to +\infty}{\rightarrow} 0.
\end{aligned}
\end{equation*}
Finally, we have
\begin{equation*}
\int_{Q_T}\dot\zeta(t)\theta^{\epsilon_k}(x-\zeta(t))u^{\epsilon_k}_x(t,x)v_t(t,x)\dd t\dd x
\underset{k\to +\infty}{\rightarrow}
\int_0^T\dot\zeta(t) u_x(t,\zeta(t))v_t(t,\zeta(t))\dd t.
\end{equation*}
\eproof

\begin{remark}
The right hand side of \eqref{AtHePe:Prop5_6} is formally equal to 
\begin{equation}
\label{AtHePe:16}
\begin{aligned}
&-u_0(\zeta(0))v_t(0,\zeta(0))-\int_0^Tu(t,\zeta(t))v_{tt}(t,\zeta(t))\dd t\dd x
\\&-\int_0^T\dot\zeta(t) u_x(t,\zeta(t))v_t(t,\zeta(t))\dd t
-\int_0^T\dot\zeta(t)u(t,\zeta(t)v_{xt}(t,\zeta(t))\dd t
\\&=-u_0(\zeta(0))v_t(0,\zeta(0))-\int_0^T\frac{\dd }{\dd t}(u(t,\zeta(t))v_t(t,\zeta(t)))\dd t
\\&+\int_0^T u_t(t,\zeta(t))v_t(t,\zeta(t))\dd t
=\int_0^T u_t(t,\zeta(t)) v_t(t,\zeta(t))\dd t.
\end{aligned}
\end{equation}
Nevertheless, with regularity of Corollary \ref{AtHePe:Cor2}, the trace  
\(\int_0^T u_t(t,\zeta(t)) v_t(t,\zeta(t))\dd t\) is not well defined.
\end{remark}

\bbproof {\emph{of Theorem \ref{AtHePe:thm1}}}.
The proof of \eqref{AtHePe:FV} follows from \eqref{AtHePe:10faible}
and Proposition \ref{AtHePe:Prop5}. Details are omitted. It remains to show that \(u(0,\cdot)=u_0\).
To this aim, notice that by \eqref{AtHePe:conv1}, \(u^{\epsilon_k}(0,\cdot)\underset{k\to +\infty}{\rightarrow} u(0,\cdot)\)
in \(\L^2_{\mathrm{loc}}(\Er)\), and that \(u^{\epsilon_k}(0,\cdot)=u_0\star\theta^{\epsilon_k}\underset{k\to +\infty}{\rightarrow} u_0\)
in \(\L^2(\Er)\). This proves the result.
\eeproof

\section*{Acknowledgment}
O.A. was funded by the Czech Academy of Sciences within Grant RVO: 67985840 and by the "EU Strategy for the Danube Region" Grant: 8X23001.
A.H. and A.P. were supported by the Ministry for Europe and Foreign Affairs (France) within Programme Danube 2023-2025 (Project PHC Danube \(n^{\rm o}\)49916UM).

\end{document}